\documentclass[journal]{IEEEtran}
\usepackage{graphicx}
\usepackage{epsfig}
\usepackage{amsmath}
\usepackage{amssymb}
\usepackage{subfigure}
\usepackage{wrapfig}
\usepackage{times,color}
\usepackage{cite,footnote,xspace,syntonly,algorithm,algorithmic,bm, url}
\input{mysymbol.sty}

\ifCLASSINFOpdf
\else
\fi
\newcommand{\calN}{{\cal N}}

\newcommand{\calJ}{{\cal J}}

\def \cN {\mathcal{N}}
\def \cL {\mathcal{L}}

\def \cP {\mathcal{P}}
\def \cJ {\mathcal{J}}

\def \bY {{\bf Y}}
\def \bS {{\bf S}}

\def \bX {{\bf X}}
\def \bA {{\bf A}}
\def \bV {{\bf V}}
\def \bW {{\bf W}}
\def \bF {{\bf F}}

\def \bI {{\bf I}}

\def \bM {{\bf M}}

\def \bV {{\bf V}}
\def \bE {{\bf E}}

\def \bZ {{\bf Z}}
\def \bW {{\bf W}}

\def \bL {{\bf P}}
\def \bQ {{\bf Q}}

\def \bO {{\bf O}}

\def \bx {{\bf x}}

\def \bl {{\bf p}}
\def \bq {{\bf q}}

\def \vec {\text{vec}}
\def \unvec {\text{unvec}}

\begin{document}
%
\title{Load Curve Data Cleansing and Imputation\\ via Sparsity and Low Rank}
%
%
%

\author{Gonzalo~Mateos,~\IEEEmembership{Member,~IEEE,}
        and~Georgios~B.~Giannakis,~\IEEEmembership{Fellow,~IEEE}
\thanks{Part of this paper will
be presented at the {\it Proc. of 3rd Intl. Conf. on Smart Grid Communications}, 
Tainan, Taiwan, Nov. 5-8, 2012.}
\thanks{The authors are with the Dept. of ECE and the Digital Technology Center, 
University of Minnesota, Minneapolis,
MN, 55455 USA (e-mail: mate0058@umn.edu; georgios@umn.edu).}}

%
%

\markboth{IEEE Transactions on Smart Grid,~Vol.~X, No.~X, XXXXXX~2012}%
{Mateos and Giannakis: Load Curve Data Cleansing and Imputation via Sparsity and Low Rank}
%



\maketitle

\begin{abstract}
The smart grid vision is to build an intelligent power network with an
unprecedented level of situational awareness and controllability over its services
and infrastructure. This paper advocates statistical inference methods to 
robustify power monitoring tasks against the outlier effects owing to faulty readings
and malicious attacks, as well as against missing data due to privacy concerns and 
communication errors. In this context, a novel 
{\em load cleansing and imputation} scheme is developed 
leveraging the low intrinsic-dimensionality of spatiotemporal load profiles and
the sparse nature of ``bad data.'' A robust estimator based on
principal components pursuit (PCP) is adopted, which effects a twofold sparsity-promoting 
regularization through an $\ell_1$-norm of the outliers, and the nuclear
norm of the nominal load profiles.
Upon recasting the non-separable nuclear norm into a form amenable
to decentralized optimization, a \textit{distributed} (D-) PCP algorithm is developed to carry out
the imputation and cleansing tasks using networked devices
comprising the so-termed advanced metering infrastructure. If D-PCP converges 
and a qualification inequality is satisfied, the novel distributed estimator provably
attains the performance of its centralized PCP counterpart, which has access to 
all networkwide data. Computer simulations and tests with real load curve data 
corroborate the convergence and effectiveness of the novel D-PCP
algorithm.
\end{abstract}


%
\IEEEpeerreviewmaketitle


\section{Introduction}
\label{sec:intro}

The US power grid has been recognized as
the most important engineering achievement of the 20th
century~\cite{nae-report}, yet it presently faces major
challenges related to efficiency, reliability, security, environmental impact,
sustainability, and market diversity issues~\cite{doe-intro}.
The crystallizing vision of the smart grid (SG) aspires to build a cyber-physical network that
can address these challenges by capitalizing on state-of-the-art 
information technologies in sensing, control, 
communication, optimization, and machine learning.
Significant effort and investment are being committed 
to architect the necessary infrastructure by installing
advanced metering systems, and establishing data communication
networks throughout the grid. Accordingly, algorithms that optimally exploit the
pervasive sensing and control capabilities of the envisioned advanced metering infrastructure
(AMI) are needed to make the necessary 
breakthroughs in the key problems in power grid monitoring and energy management.
This is no  easy endeavor though, in view of the challenges posed by 
increasingly distributed network operations under strict reliability 
requirements, also facing  malicious cyber-attacks.

{\em Statistical inference} techniques are expected to play an increasingly
instrumental role in power system monitoring~\cite{GE_spmag}, not only to meet the anticipated 
``big data" deluge as the installed base of phasor measurement units (PMUs) 
reaches out throughout the grid, but also to robustify the monitoring tasks against the
``outlier" effects owing to faulty readings, malicious attacks, and
communication errors, as well as against missing data due to privacy concerns and 
technical anomalies~\cite{Scharf}. In this context, a 
{\em load cleansing and imputation} scheme is developed in this paper, building on recent 
advances in sparsity-cognizant information processing~\cite{elements_of_statistics}, 
low-rank matrix completion~\cite{F02}, and large-scale distributed optimization~\cite{Boyd_ADMM}.

\textit{Load curve} data refers to the electric energy consumption periodically recorded
by meters at points of interest across the power grid, e.g., end-user premises, buses, and substations.
Accurate load profiles are critical assets aiding operational decisions
in the envisioned SG system~\cite{cllce10tsg}, and are essential for
load forecasting~\cite{Shah-MO}.
However, in the process of acquiring and transmitting such
massive volumes of information for centralized processing, data are oftentimes corrupted or
lost altogether. In a smart monitoring context for instance, incomplete load profiles emerge
due to three reasons: (r1) PMU-instrumented buses are few; (r2) SCADA data become
available at a considerably smaller time scale
than PMU data; and (r3) regional operators are not
willing to share all their variables \cite{rso11psce}. Moreover, a major
requirement for grid monitoring is robustness to outliers,
i.e., data not adhering to nominal models~\cite{Abur_book,mili}. Sources of
so-termed ``bad data'' include meter failures, as well as strikes,
unscheduled generator shutdowns, and extreme weather
conditions~\cite{cllce10tsg,glliw10tps}. Inconsistent data can also be due to malicious
(cyber-) attacks that induce abrupt load changes, 
or counterfeit meter readings~\cite{lnr09cccf}. 

In light of the aforementioned observations, the first contribution of this paper 
is on modeling spatiotemporal load profiles, accounting for the 
structure present to effectively impute missing data and 
devise robust load curve estimators stemming from 
convex optimization criteria (Section \ref{sec:prob_state}). 
Existing approaches to load curve cleansing have relied on separate processing
per time series~\cite{cllce10tsg,gm_gg_2012,glliw10tps}, 
and have not capitalized on spatial correlations to improve performance. 
The aim is for \textit{minimal-rank} cleansed load data, while also exploiting outlier 
\textit{sparsity} across buses and time. An estimator tailored
to these specifications is principal components pursuit (PCP)~\cite{CLMW09,CSPW11,zlwcm10},
which is outlined in Section \ref{sec:PCP}. PCP minimizes a tradeoff
between the least-squares (LS) model fitting error and a twofold
sparsity-promoting regularization, implemented through an $\ell_1$-norm of the outliers 
and the nuclear norm of the nominal load profiles. While PCP has been widely adopted
in computer vision~\cite{CLMW09}, for voice separation in music~\cite{hcsh12}, and unveiling 
network anomalies~\cite{tsp_rankminimization_2012}, its benefits to power
systems engineering and monitoring remains so far largely unexplored. 
The second contribution pertains to developing a \textit{distributed} (D-) PCP algorithm, to carry out
the imputation and cleansing task using a network of interconnected devices
as part of the AMI (Section \ref{sec:dist}). This is possible by leveraging a general
algorithmic framework for sparsity-regularized rank minimization put forth in
\cite{tsp_rankminimization_2012}. Upon recasting the (non-separable) nuclear norm
present in the PCP cost into a form amenable to decentralized optimization (Section
\ref{ssec:separable}), the D-PCP iterations are obtained in 
Section \ref{ssec:ad_mom} via the multi-block alternating-directions method of multipliers (ADMM)
solver~\cite{Bertsekas_Book_Distr,Boyd_ADMM}. In a nutshell, per iteration 
each smart meter exchanges simple messages with its (directly connected)
neighbors in the network, and then solves its own optimization problem to refine
its current estimate of the cleansed load profile. In the context of power systems, the ADMM
has been recently adopted to carry out dynamic network energy management~\cite{boyd_dist_energy}, and 
distributed robust state estimation~\cite{kekatos}. 

Computer simulations corroborate the convergence and optimality of the novel D-PCP
algorithm, and demonstrate its effectiveness in cleansing and imputing 
real load curve data (Section \ref{sec:sims}). Concluding remarks 
and directions for future research are outlined in Section \ref{sec:conc}, while
a few algorithmic details are deferred to the Appendix.


\section{Modeling and Problem Statement}
\label{sec:prob_state}

This section introduces the model for
(possibly) incomplete and grossly corrupted load curve measurements,
acquired by geographically-distributed metering devices monitoring the power grid. 
The communication network model needed to account for exchanges of information
among smart meters is described as well. Lastly, 
the task of load curve cleansing and imputation is formally stated.


\subsection{Spatiotemporal load curve data model}
\label{ssec:data_model}

Let the $N\times 1$ vector $\mathbf{y}(t):=[y_{1,t},\ldots,y_{N,t}]'$ ($'$ stands
for transposition)
collect the spatial load profiles measured by smart meters monitoring $N$ network
nodes (buses, residential premises), at a given discrete-time instant  $t\in[1,T]$. Consider 
the $N\times T$ matrix of observations 
$\mathbf{Y}:=[\mathbf{y}(1),\ldots,\mathbf{y}({ T})]$. The $n$-th row $(\mathbf{y}_n)'$ of $\mathbf{Y}$
is the time series of energy consumption (load curve) measurements at node $n$, while
the $t$-th column $\mathbf{y}(t)$ of $\bY$ represents a snapshot of the networkwide loads
taken at time $t$.  To model missing data, consider the set
$\Omega\subseteq \{1,\ldots,N\}\times\{1,\ldots,{ T}\}$ of index pairs
$(n,t)$ defining a sampling of the entries of $\mathbf{Y}$. Introducing
the matrix sampling operator $\mathcal{P}_\Omega(\cdot)$, which sets the entries 
of its matrix argument not indexed by $\Omega$ to zero and leaves the rest
unchanged, the (possibly) incomplete spatiotemporal load curve data in the presence of outliers
can be modeled as
\begin{equation}
\label{eq:data_model}
\mathcal{P}_\Omega(\mathbf{Y}) = \mathcal{P}_\Omega(\mathbf{X} + \mathbf{O} + \mathbf{E})
\end{equation}
where $\mathbf{X}$, $\mathbf{O}$, and $\mathbf{E}$ denote
the nominal load profiles, the outliers, and small measurement errors, respectively.
For nominal observations $y_{n,t}=x_{n,t}+e_{n,t}$, one has $o_{n,t}=0$.

\begin{remark}[Model under-determinacy]\label{rem:underdetermined}
\normalfont The model is inherently under-determined, since even for the (most favorable) 
case of full data, i.e., $\Omega\equiv\{1,\ldots,N\}\times\{1,\ldots,{ T}\}$, there are twice as
many unknowns in $\mathbf{X}$ and $\mathbf{O}$ as there is data in $\mathbf{Y}$. 
Estimating $\mathbf{X}$ and $\mathbf{O}$ becomes even more challenging
when data are missing, since the number of unknowns remains the same, but the amount
of data is reduced. 
\end{remark}

In any case, estimation of $\{\mathbf{X},\mathbf{O}\}$
from $\mathcal{P}_\Omega(\mathbf{Y})$ is an ill-posed problem unless one
introduces extra structural assumptions on the model components to reduce
the effective degrees of freedom. To this end, two 
cardinal properties of $\mathbf{X}$ and $\mathbf{O}$ will prove instrumental.
First, common temporal patterns among the energy consumption of
a few broad classes of loads (e.g., industrial, residential, seasonal)
in addition to their (almost) periodic behaviors render most rows
and columns of $\mathbf{X}$ linearly dependent, and thus $\mathbf{X}$ typically has
\emph{low-rank}. Second, outliers (or attacks) only occur sporadically in
time and affect only a few buses, yielding a \emph{sparse} matrix
$\mathbf{O}$. Smoothness of the nominal load curves is related
to the low-rank property of $\mathbf{X}$, which was adopted in~\cite{cllce10tsg} to 
motivate a smoothing splines-based algorithm for cleansing. However, while in~\cite{cllce10tsg}
smoothness was enforced per load profile time-series, i.e.,
per row $(\bx_n)'$ of $\bX$; the low-rank property of $\mathbf{X}$ also
captures the spatial dependencies introduced by the network.  Approaches 
capitalizing on outlier- and ``bad data-'' sparsity can be found in e.g.,~\cite{tong_bad_data,kekatos} 
and~\cite{gm_gg_2012}.


\subsection{Communication network model}
\label{ssec:comm_model}

Suppose that on top of the energy measurement functionality, 
the $N$ networked smart meters are capable of performing simple local computations, as well as 
exchanging messages among directly connected neighbors. Single-hop
communication models are appealing due to their simplicity, 
since one does not have to incur the routing overhead.
The AMI network is naturally abstracted to an
undirected graph $G(\cN,\cL)$, where the vertex set
$\cN:=\{1,\ldots,N\}$ corresponds to the network nodes, and the
edges (links) in $\cL$ represent pairs of nodes that are
connected via a physical communication channel. 
Node $n\in\cN$ communicates with its single-hop
neighboring peers in $\cJ_n$, and the size of the neighborhood will
be henceforth denoted by $|\cJ_n|$. The graph $G$ is assumed
connected, i.e., there exists a (possibly multihop) path that joins
any pair of nodes in the network. This requirement ensures that
the network is devoid of multiple isolated (connected) components, and allows
for the data collected by e.g., smart meter $n$, namely the $n$-th row $(\mathbf{y}_n)'$
of $\mathbf{Y}$, to eventually reach every other node in the network. This
way, even when only local interactions are allowed, the flow of information can percolate
the network. 

The importance of the network model will become apparent in Section \ref{sec:dist}.


\subsection{Load curve cleansing and imputation}
\label{ssec:cleanse_impute}

The load curve cleansing and imputation problem studied here
entails identification and removal of outliers (or ``bad data''), in addition
to completion of the missing entries from the nominal load matrix, and denoising of the observed ones.
To some extent, it is a joint estimation-interpolation (prediction)-detection problem.
With reference to \eqref{eq:data_model}, given
generally incomplete, noisy and outlier-contaminated spatiotemporal
load data $\mathcal{P}_{\Omega}(\mathbf{Y})$, the cleansing and imputation tasks amount to 
estimating the nominal load profiles $\mathbf{X}$ and the outliers $\mathbf{O}$,
by leveraging the low-rank property of $\mathbf{X}$ and the sparsity in $\mathbf{O}$.
Collaboration between metering devices (collecting networkwide data) is considered here,
rather than local processing of load curves per bus.

Note that load cleansing and imputation are different from \textit{load forecasting}~\cite{Shah-MO},
which amounts to predicting future load demand based on historical data
of energy consumption and the weather conditions. Actually, cleansing and imputation
are critical preprocessing tasks utilized to enhance the quality of load data, that
would be subsequently used for load forecasting and optimum power flow~\cite{PSA_book}.


\section{Prinicipal Components Pursuit}
\label{sec:PCP}

An estimator matching nicely the specifications of the load
curve cleansing and imputation problem stated in Section \ref{ssec:cleanse_impute}, is the so-termed  
(stable) principal components pursuit (PCP)~\cite{CLMW09,CSPW11,zlwcm10},
that will be outlined here for completeness. 
PCP seeks estimates $\{\hat{\mathbf{X}}, \hat{\mathbf{O}}\}$ as the minimizers of
\begin{align}
\text{(P1)}~~~~~\min_{\{\mathbf{X},
\mathbf{O}\}} \frac{1}{2}\left\|\mathcal{P}_\Omega(\mathbf{Y} - \mathbf{X}-\mathbf{O}
)\right\|_F^2 + \lambda_* \left\| \mathbf{X}\right\|_{*}+ \lambda_1 \left\|\mathbf{O}\right\|_1
\nonumber
\end{align}
where the $\ell_1$-norm $\|\bO\|_1:=\sum_{n,t}|o_{n,t}|$ and the nuclear norm
$\|\bX\|_*:=\sum_{i}\sigma_{i}(\bX)$
($\sigma_{i}(\bX)$ denotes the $i$-th singular value of $\bX$) are utilized to promote sparsity in the
number of outliers (nonzero entries) in $\mathbf{O}$, and the low rank
of $\mathbf{X}$, respectively. The nuclear and $\ell_1$-norms
are the closest convex surrogates to the rank and cardinality functions,
which albeit the most natural criteria they are in general NP-hard to 
optimize. The tuning parameters 
$\lambda_1,\lambda_*\geq 0$ control the tradeoff
between fitting error, rank, and sparsity level of the solution. 
When an estimate $\hat{\sigma}_v^2$ of the observation noise variance is available,
guidelines for selecting $\lambda_*$ and $\lambda_1$ have been proposed in~\cite{zlwcm10}. 
The nonzero entries in $\hat{\mathbf{O}}$ reveal ``bad data'' across both buses and time.
Clearly, it does not make sense to flag outliers in data that
has not been observed, namely for $(n,t)\notin\Omega$.
In those cases (P1) yields $\hat{o}_{n,t}=0$ since both the Frobenious and 
$\ell_1$-norms are separable across the entries of their matrix arguments.

Being convex (P1) is computationally appealing, and it has been shown to attain good performance 
in theory and practice. For instance, in the absence of noise and when there is 
no missing data, identifiability and exact recovery conditions 
were reported in~\cite{CLMW09} and~\cite{CSPW11}.
Even when data are missing, it is possible to recover the low-rank component
under some technical assumptions~\cite{CLMW09}. Theoretical 
performance guarantees in the presence of noise are also available~\cite{zlwcm10}.

Regarding algorithms, a PCP solver based on the accelerated proximal gradient
method was put forth in~\cite{rpca_proximal}, while the ADMM
was employed in~\cite{rpca_admom}. 
Implementing these \textit{centralized} algorithms presumes that networked metering devices 
continuously communicate their local load measurements $y_{n,t}$ to a 
central monitoring and data analytics station, which uses their aggregation in 
$\mathcal{P}_\Omega(\mathbf{Y})$ to reject outliers and impute missing data.
While for the most part this is the prevailing operational paradigm nowadays, 
there are limitations associated with this architecture.
For instance, collecting all this information centrally may lead to excessive overhead
in the communication network, especially
when the rate of data acquisition is high at the meters. Moreover, minimizing 
(or avoiding altogether) the exchanges of raw
measurements may be desirable for privacy and cyber-security reasons, as well as 
to reduce unavoidable communication errors that translate to missing
data. Performing the optimization in a centralized fashion raises robustness concerns as well, since the
central data analytics station represents an isolated point of failure.
These reasons motivate devising fully-distributed iterative algorithms for PCP, 
embedding the load cleansing and imputation functionality to the AMI. This is the 
subject of the next section.


\section{Distributed Cleansing and Imputation}
\label{sec:dist}

A distributed (D-)PCP algorithm to solve (P1) using a network of smart meters
(modeled as in Section \ref{ssec:comm_model}) should be
understood as an iterative method, whereby each node carries out simple local (optimization)
tasks per iteration $k=1,2,\ldots$, and exchanges messages only with its directly 
connected neighbors. The ultimate goal is for each node to form local estimates
$\mathbf{x}_n[k]$ and $\bbo_n[k]$ that coincide with the $n$-th rows of $\hat{\mathbf{X}}$
and $\hat{\mathbf{O}}$ as $k\to\infty$, where $\{\hat{\mathbf{X}},\hat{\mathbf{O}}\}$
is the solution of (P1) obtained when all data $\mathcal{P}_{\Omega}(\mathbf{Y})$
are centrally available. Attaining the centralized performance with distributed data is
impossible if the network is disconnected.

To facilitate reducing the computational complexity and memory
storage requirements of the D-PCP algorithm sought, it is
henceforth assumed that an upper bound $\textrm{rank}(\hat\bX)\leq
\rho$ is a priori available [recall $\hat\bX$ is the estimated low-rank
cleansed load profile obtained via (P1)]. As argued next, the smaller the value of
$\rho$, the more efficient the algorithm becomes. 
Small values of $\rho$ are well motivated due to the low intrinsic 
dimensionality of the spatiotemporal load profiles (cf. Section \ref{ssec:data_model}). 
Because $\textrm{rank}(\hat\bX)\leq \rho$, (P1)'s search space is
effectively reduced and one can factorize the decision variable as
$\bX=\bL\bQ'$, where $\bL$ and $\bQ$ are $N \times \rho$ and $T
\times \rho$ matrices, respectively. Adopting this reparametrization
of $\bX$ in (P1) and making explicit the distributed nature of the
data (cf. Section \ref{ssec:comm_model}), one arrives at an 
equivalent optimization problem
\begin{align}
\text{(P2)}~~~~~\min_{\{\bL,\bQ,\bO\}}&~\sum_{n=1}^N\left[\frac{1}{2}\|\cP_{\Omega_n}(\mathbf{y}_n 
-\bQ\bl_n-\mathbf{o}_n)\|_2^2
\right.\nonumber\\
&\hspace{0.8cm}\left. +
\frac{\lambda_*}{N}\|\bL\bQ'\|_{\ast}+\lambda_1\|\mathbf{o}_n\|_1\right]\nonumber
\end{align}
which is non-convex due to the bilinear term $\bL\bQ'$, and where
$\bL :=\left[\bl_1,\ldots,\bl_N\right]^\prime$. The number of variables is 
reduced from $2NT$ in (P1), to $\rho(N + T) + NT$ in (P2). 
The savings can be significant when $\rho$ is small, and both $N$ and $T$
are large. Note that the dominant $NT$-term in the variable count of 
(P2) is due to $\mathbf{O}$, which is sparse and
can be efficiently handled even when both $N$ and $T$ are large. 

\begin{remark}[Challenges facing distributed implementation]\label{remark:challenges}
\normalfont Problem (P2) is still not amenable for distributed implementation
due to: (c1) the non-separable nuclear norm present in the cost
function; and (c2) the global variable $\bQ$ coupling the per-node
summands. 
\end{remark}

Challenges (c1)-(c2) are dealt with in the ensuing sections.


\subsection{A separable low-rank regularization}
\label{ssec:separable}

To address (c1), consider the following alternative characterization of
the nuclear norm (see e.g.~\cite{tsp_rankminimization_2012})
\begin{equation}\label{eq:nuc_nrom_def}
\|\bX\|_*:=\min_{\{\bL,\bQ\}}~~~ \frac{1}{2}\left(\|\bL\|_F^2+\|\bQ\|_F^2 \right),\quad
\text{s. to}~~~ \bX=\bL\bQ'.
\end{equation}
The optimization \eqref{eq:nuc_nrom_def} is over all possible bilinear factorizations
of $\bX$, so that the number of columns $\rho$ of $\bL$ and
$\mathbf{Q}$ is also a variable. Leveraging \eqref{eq:nuc_nrom_def}, 
the following reformulation of (P2) provides an important first step towards obtaining 
the D-PCP algorithm:
\begin{align}
\text{(P3)}~~~\min_{\{\bL,\bQ,\bO\}}& \sum_{n=1}^N \left[ \frac{1}{2}\|\cP_{\Omega_n}(\mathbf{y}_n 
-\bQ\bl_n-\mathbf{o}_n)\|_2^2+ \lambda_1\|\mathbf{o}_n\|_1\right.\nonumber\\
&\hspace{0.8cm}\left.+
\frac{\lambda_{*}}{2N}\left(N\|\bl_n\|_2^2 + \|\bQ\|_F^2 \right) \right].\nonumber
\end{align}
As asserted in~\cite[Lemma 1]{tsp_rankminimization_2012}, adopting the separable Frobenius-norm 
regularization in (P3) comes with no loss of optimality relative to (P1), 
provided $\textrm{rank}(\hat\bX)\leq\rho$. By finding the global minimum of (P3)
[which could have considerably less variables than (P1)], 
one can recover the optimal solution of (P1). This could be challenging however, since (P3) is
non-convex and it may have stationary points which need not be globally
optimum. 

Interestingly, it is possible to certify global optimality 
of a stationary point $\{\bar{\bL},\bar{\bQ},\bar{\bO}\}$ of (P3). Specifically, one
can establish that if
$\|\cP_{\Omega}(\bY - \bar{\bL}\bar{\bQ}' - \bar{\bO})\| < \lambda_*$, then
$\{\hat\bX:=\bar{\bL}\bar{\bQ}',\hat{\bO}:=\bar{\bO}\}$ is the globally optimal solution of
(P1)~\cite[Prop. 1]{tsp_rankminimization_2012}.
The qualification condition $\|\cP_{\Omega}(\bY - \bar{\bL}\bar{\bQ}' - \bar{\bO})\| < \lambda_*$
captures tacitly the role of $\rho$. In particular, for sufficiently small $\rho$ the residual
$\|\cP_{\Omega}(\bY - \bar{\bL}\bar{\bQ}' - \bar{\bO})\|$
becomes large and consequently the condition is violated [unless $\lambda_*$ is large
enough, in which case a sufficiently low-rank solution to (P1) is expected]. The condition
on the residual also implicitly enforces $\textrm{rank}(\hat\bX) \leq \rho$, which is necessary for
the equivalence between (P1) and (P3).


\subsection{Local variables and consensus constraints}
\label{ssec:constraints}

To decompose the cost in (P3), in which summands inside the square brackets 
are coupled through the global variable $\bQ$ [cf. (c2) under Remark 
\ref{remark:challenges}], introduce auxiliary variables $\{\bbQ_n\}_{n=1}^N$
representing local estimates of $\bQ$ per smart meter $n$. To obtain
a separable PCP formulation, use these
estimates along with \textit{consensus} constraints
\begin{align}
\text{(P4)}~~~~\min_{\{\bL_n,\bQ_n,\bO\}}& \sum_{n=1}^N \left[\frac{1}{2}\|\cP_{\Omega_n}(\mathbf{y}_n 
-\bQ_n\bl_n-\mathbf{o}_n)\|_2^2\right.\nonumber\\
&\left.\hspace{0.1cm}  + \lambda_1\|\mathbf{o}_n\|_1+
\frac{\lambda_{*}}{2N}\left(N\|\bl_n\|_2^2 + \|\bQ_n\|_F^2 \right) \right] \nonumber\\
\text{s. to} &\quad \bQ_n=\bQ_m,\quad m\in\cJ_n,\:n\in\mathcal{N}.\nonumber
\end{align}
Notice that (P3) and (P4) are equivalent optimization problems,
since the network graph $G(\cN,\cL)$ is connected by assumption. 
Even though consensus is a fortiori imposed only within neighborhoods, it extends to
the whole (connected) network and local estimates agree on the
global solution of (P3). To arrive at the desired D-PCP
algorithm, it is convenient to reparametrize the consensus
constraints in (P4) as
\begin{align}
\bQ_n={}&{}\bar{\bF}_n^m,\:\bQ_m=\tilde{\bF}_n^m, \textrm{ and }
\bar{\bF}_n^m=\tilde{\bF}_n^m, \quad m\in\cJ_n,\:n\in\mathcal{N}\label{eq:constr_1}
\end{align}
where $\{\bar{\bF}_n^m,\tilde{\bF}_n^m\}_{n\in\cN}$, are auxiliary
optimization variables that will be eventually eliminated (cf. Remark \ref{rem:simplific}).


\subsection{The D-PCP algorithm}
\label{ssec:ad_mom}

To tackle (P4), associate Lagrange multipliers $\bar{\bM}_n^m$ and
$\tilde{\bM}_n^m$ with the first pair of consensus constraints in
\eqref{eq:constr_1}. Introduce the quadratically \textit{augmented}
Lagrangian function~\cite{Bertsekas_Book_Distr}
\begin{align}\label{augLagr}
\ccalL_c\left(\mathcal{V}_1,\mathcal{V}_2,\mathcal{V}_3,\mathcal{M}\right){}={}&
\sum_{n=1}^N\left[\frac{1}{2}\|\cP_{\Omega_n}(\mathbf{y}_n 
-\bQ_n\bl_n-\mathbf{o}_n)\|_2^2\right.\nonumber\\
&\hspace{-3cm}\left. + \lambda_1\|\mathbf{o}_n\|_1+\frac{\lambda_{*}}{2N}\left(N\|\bl_n\|_2^2 + \|\bQ_n\|_F^2 
\right)\right]\nonumber\\
&\hspace{-3cm} + \sum_{n=1}^N\sum_{m\in\calJ_n}\left(\langle \bar{\bM}_n^m,\bQ_n-\bar{\bF}_n^m\rangle
+\langle \tilde{\bM}_n^m,\bQ_m-\tilde{\bF}_n^m\rangle \right)\nonumber\\
&\hspace{-3cm}+\frac{c}{2}\sum_{n=1}^N\sum_{m\in\calJ_n}\left(\|\bQ_n-\bar{\bF}_n^m\|_F^2
+\|\bQ_m-\tilde{\bF}_n^m\|_F^2\|_F^2\right)
\end{align}
where $c>0$ is a penalty parameter, and the primal variables are split
into three groups $\mathcal{V}_1:=\{\bQ_n\}_{n=1}^{N}$, $\mathcal{V}_2:=\{\bl_n\}_{n=1}^{N}$ and
$\mathcal{V}_3:=\{\bbo_n,\bar{\bF}_n^m$ $,
\tilde{\bF}_n^m\}_{n\in\cN}^{m\in\cJ_n}$. For notational brevity,
collect all Lagrange multipliers in
$\mathcal{M}:=\{\bar{\bM}_n^m,\tilde{\bM}_n^m\}_{n\in\cN}^{m\in\cJ_n}$.
Note that the remaining constraints in \eqref{eq:constr_1}, namely $
\mathcal{C}_F:=\{\bar{\bF}_n^m=\tilde{\bF}_n^m,
\:m\in\cJ_n,\:n\in\mathcal{N}\}$, have not been dualized.

To minimize (P4) in a distributed fashion, (a multi-block variant of) the ADMM 
will be adopted here. The ADMM is
an iterative augmented Lagrangian method especially well suited for
parallel processing~\cite{Bertsekas_Book_Distr,Boyd_ADMM}, which has been
proven successful to tackle the optimization tasks stemming from
general distributed estimators of deterministic and (non-)stationary
random signals; see e.g.,~\cite{kekatos,tsp_rankminimization_2012}
and references therein.
The proposed solver entails an iterative procedure comprising four steps per
iteration $k=1,2,\ldots$, 
[$n\in\cN,\:m\in\cJ_n$ in
\eqref{eq:multi_barC}-\eqref{eq:multi_tildeC}]
\begin{description}
\item [{\bf [S1]}] \textbf{Update dual variables:}\begin{align}
  \bar{\bM}_n^m[k]&=\bar{\bM}_n^m[k-1]+c (\bQ_n[k]-\bar{\bF}_n^m[k])\label{eq:multi_barC}\\
\tilde{\bM}_n^m[k]&=\tilde{\bM}_n^m[k-1]+c (\bQ_m[k]-\tilde{\bF}_n^m[k]).\label{eq:multi_tildeC}
\end{align}

\item [{\bf [S2]}]  \textbf{Update first group of primal variables:}
    \begin{equation}\hspace{-1cm}\mathcal{V}_1[k+1]=
    \mbox{arg}\:\min_{\mathcal{V}_1}\ccalL_c\left(\mathcal{V}_1,\mathcal{V}_2[k],\mathcal{V}_3[k],
    \mathcal{M}[k]\right).
    \label{S2_ADMOM}\end{equation}

\item [{\bf [S3]}]  \textbf{Update second group of primal variables:}
        \begin{equation}\hspace{-1cm}\mathcal{V}_2[k+1]=
            \mbox{arg}\:\min_{\mathcal{V}_2}\ccalL_c\left(\mathcal{V}_1[k+1],\mathcal{V}_2,\mathcal{V}_3[k],
            \mathcal{M}[k]\right).
        \label{S3_ADMOM}\end{equation}

\item [{\bf [S4]}]  \textbf{Update third group of primal variables:}
        \begin{equation}\hspace{-1cm}\mathcal{V}_3[k+1]=
            \mbox{arg}\:\min_{\mathcal{V}_3\in \mathcal{C}_F}\ccalL_c\left(\mathcal{V}_1[k+1],\mathcal{V}_2[k+1]
            ,\mathcal{V}_3,\mathcal{M}[k]\right)
        \label{S4_ADMOM}\end{equation}
\end{description}
which amount to a block-coordinate descent method cycling over $\mathcal{V}_1\rightarrow
\mathcal{V}_2\rightarrow\mathcal{V}_3$ to minimize $\ccalL_c$,
and dual variable updates~\cite{Bertsekas_Book_Distr}. At each step while minimizing the
augmented Lagrangian, the variable groups not being updated are treated as
fixed, and are substituted with their most up to date values. Different
from the standard two-block ADMM~\cite{Bertsekas_Book_Distr,Boyd_ADMM}, 
the multi-block variant here cycles over three groups of primal 
variables~\cite{luo_ADMM}.

Reformulating the estimator (P1) to its equivalent form
(P4) renders the augmented Lagrangian in \eqref{augLagr}
highly decomposable. The separability comes in two flavors,
both with respect to the variable groups $\mathcal{V}_1$-$\mathcal{V}_3$, as well as across
the network nodes $n\in\cN$. This leads to highly
parallelized, simplified recursions to be run by the networked
smart meters. Specifically, it is shown in the Appendix
that the aforementioned ADMM steps [S1]-[S4] give rise to the D-PCP iterations tabulated
under Algorithm \ref{tab:table_2}. Per iteration, each device updates: [S1] a local
matrix of dual prices $\bS_n[k]$;
[S2]-[S3] local cleansed load estimates $\bQ_n[k+1]$ and $\bl_n[k+1]$ obtained
as solutions to respective unconstrained quadratic problems (QPs); and [S4] its local outlier
vector, through a sparsity-promoting soft-thresholding operation. 
The $(k+1)$-st iteration is concluded after smart meter $n$ transmits $\bQ_n[k+1]$
to its single-hop neighbors in $\mathcal{J}_n$. Regarding communication cost,
$\bQ_n[k+1]$ is a $T \times \rho$ matrix and its transmission does not
incur significant overhead for small $\rho$. Observe also that
$\mathcal{P}_{\Omega}(\mathbf{y}_n)$ need not be exchanged which is 
desirable to preserve data secrecy, and the communication cost
is independent of $N$.

\begin{algorithm}[t]
\caption{: D-PCP at smart meter $n\in\cN$} \small{
\begin{algorithmic}
    \STATE \textbf{input} $\bby_n, \bm{\Omega}_n, \lambda_{*}, \lambda_1,$ and $c$.
    \STATE \textbf{initialize}
$\bS[0]=\mathbf{0}_{T\times \rho}$, and $\bQ_n[1],\:\bl_n[1]$ at random.
    \FOR {$k=1,2$,$\ldots$}
        \STATE Receive $\{\bQ_m[k]\}$ from neighbors $m\in\cJ_n$.
        \STATE \textbf{[S1] Update local dual variables:}
        \STATE $\bS_n[k]=\bS_n[k-1] + c \sum_{m\in\cJ_n}(\bQ_n[k]-\bQ_m[k])$.
        \STATE \textbf{[S2] Update first group of local primal variables:}
        \STATE $\bA_n[k+1]:=\left\{ (\bl_n[k]\bl_n'[k]) \otimes \bm{\Omega}_n +
                        (\lambda_*/N + 2 c |\cJ_n|) \bI_{\rho T}\right\}^{-1}$
        \STATE
        $\bQ_n[k+1]{}={} \unvec \Big( \bA_n[k+1]\Big\{ 
        (\bl_n[k] \otimes \bm{\Omega}_n) (\bby_n-\bbo_n[k])$
        \STATE  \hspace{1cm}$\left.\left.- \vec(\bS_n[k]) +  \vec(c\sum_{m\in\cJ_n}(\bQ_n[k] +
                                \bQ_m[k])) \right\}\right).$
        \STATE \textbf{[S3] Update second group of local primal variables:}
        \STATE
        $\bl_n[k+1] = \left\{\bQ_n'[k+1]\bm{\Omega}_n\bQ_n[k+1]+\bI_{\rho}\right\}^{-1}$
        \STATE  \hspace{1.8cm}$\times\bQ_n'[k+1]\bm{\Omega}_n(\bby_n-\bbo_n[k]).$
        \STATE \textbf{[S4] Update third group of local primal variables:}
                \STATE
                $\bbo_n[k+1] = \mathcal{S}_{\lambda_1}\left(\bm{\Omega}_n(\bby_n-
                \bQ_n[k+1]\bl_n[k+1])\right)$.
        \STATE Transmit $\bQ_n[k+1]$ to neighbors $m\in\cJ_n$.
    \ENDFOR
    \STATE \textbf{return} $\bQ_n[\infty], \bl_n[\infty], \bbo_n[\infty]$.
\end{algorithmic}}
\label{tab:table_2}
\end{algorithm}

Before moving on, a clarification on the notation used
in Algorithm \ref{tab:table_2} is due. To define matrix $\bm{\Omega}_n$
in [S2]-[S4], observe first that the 
local sampling operator can be expressed as 
$\mathcal{P}_{\Omega_n}(\mathbf{z})=\bm{\omega}_n\odot\mathbf{z}$,
where $\odot$ denotes Hadamard product, and the binary masking vector 
$\bm{\omega}_n\in\{0,1\}^T$ has entries equal to $1$ if the corresponding
entry of $\bbz$ is observed, and $0$ otherwise. It is then apparent that the Hadamard
product can be replaced with the usual matrix-vector product as 
$\mathcal{P}_{\Omega_n}(\mathbf{z})=\bm{\Omega}_n\mathbf{z}$,
where $\bm{\Omega}_n:=\textrm{diag}(\bm{\omega}_n)$. Operators $\otimes$ and $\textrm{vec}[\cdot]$
denote Kronecker product and matrix vectorization, respectively. Finally, 
the soft-thresholding operator is 
$\mathcal{S}_{\lambda_1}(\cdot):=\textrm{sign}(\cdot)\max(|\cdot|-\lambda_1,0)$.

\begin{remark}[Elimination of redundant variables]\label{rem:simplific}
\normalfont Careful inspection of Algorithm \ref{tab:table_2}
reveals that the redundant auxiliary variables
$\{\bar{\bF}_n^m,\tilde{\bF}_n^m,\tilde{\bM}_n^m\}_{n\in\cN}^{m\in\cJ_n}$
have been eliminated. Each smart meter, say the $n$-th, does not need to
\textit{separately} keep track of all its non-redundant multipliers
$\{\bar{\bM}_n^m\}_{m\in\cJ_n}$, but only update their respective
sums $\bS_n[k]:=2\sum_{m\in\cJ_n}$ $\bar{\bM}_n^m[k]$.
\end{remark}

When employed to solve non-convex problems such as (P4), ADMM
so far offers no convergence guarantees. However, there is ample
experimental evidence in the literature which supports convergence
of ADMM, especially when the non-convex problem at hand exhibits
``favorable'' structure~\cite{Boyd_ADMM}. For instance, (P4) is bi-convex and gives
rise to the strictly convex optimization subproblems
each time $\mathcal{L}_c$  is minimized with respect to one 
of the group variables, which admit unique closed-form
solutions per iteration [cf. \eqref{S2_ADMOM}-\eqref{S4_ADMOM}]. 
This observation and the linearity of the
constraints suggest good convergence properties for the D-PCP algorithm.  
Extensive numerical tests including those
presented in Section \ref{sec:sims} demonstrate that this is
indeed the case. While a formal convergence proof is the subject of
ongoing investigation, the following proposition asserts that upon convergence, the D-PCP
algorithm attains consensus and global optimality. For a proof (omitted here
due to space limitations), see~\cite[Appendix C]{tsp_rankminimization_2012}. 

\begin{proposition}\label{prop:ther_1}
Suppose iterates $\{\bQ_n[k],\bl_n[k],\bbo_n[k]\}_{n\in\cN}$ generated
by Algorithm \ref{tab:table_2} converge to
$\{\bar{\bQ}_n,\bar{\bl}_n,\bar{\bbo}_n\}_{n\in\cN}$. If $\{\hat{\bX},\hat{\bO}\}$ is the
optimal solution of (P1), then
$\bar{\bQ}_1=\bar{\bQ}_2=\ldots=\bar{\bQ}_N$. Also, if
$\|\cP_{\Omega}(\bY-\bar{\bL}{\bar{\bQ}_1}'-\bar{\bO})\| < \lambda_*$, then
$\{\hat{\bX}=\bar{\bL}\bar{\bQ}_1^\prime,\hat{\bO}=\bar{\bO}\}$.
\end{proposition}


\section{Numerical Tests}
\label{sec:sims}
This section corroborates convergence and gauges performance of the
D-PCP algorithm, when tested using synthetic and real load curve data.


\subsection{Synthetic data tests}
\label{ssec:sim_data}

\begin{figure}[t]
\centering
\includegraphics[width=\linewidth]{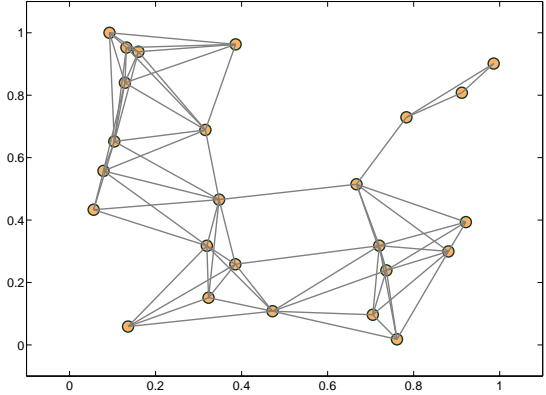}
\caption{A simulated network graph with $N=25$ nodes.}
\label{fig:Fig_1}
\vspace{-0.5cm}
\end{figure}

A network of $N=25$ smart meters is generated as a realization
of the random geometric graph model, meaning nodes are randomly placed
on the unit square and two nodes
communicate with each other if their Euclidean distance is less than a
prescribed communication range of $d_c=0.4$; see Fig. \ref{fig:Fig_1}. The time horizon is
$T=600$. Entries of $\bE$ are independent and identically distributed (i.i.d.),
zero-mean, Gaussian with variance $\sigma^2=10^{-3}$; i.e., $e_{l,t}\sim \mathcal{N}(0,\sigma^2)$.
Low-rank spatiotemporal load profiles with rank $r=3$ are
generated from the bilinear factorization model $\bX = \bW\bZ'$, where
$\bW$ and $\bZ$ are $N\times r$
and $T \times r$ matrices with i.i.d. entries drawn from Gaussian
distributions $\mathcal{N}(0,100/N)$ and $\mathcal{N}(0,100/T)$,
respectively. Every entry of $\bO$ is randomly drawn from the
set $\{-1,0,1\}$ with ${\rm Pr} (o_{n,t}=-1)={\rm Pr}(o_{n,t}=1)=5\times 10^{-2}$.
To simulate missing data, a sampling matrix $\bm{\Omega}\in\{0,1\}^{N\times T}$
is generated with i.i.d. Bernoulli distributed entries $o_{n,t}\sim\textrm{Ber}(0.7)$ 
($30$\% missing data on average). Finally, measurements are generated as 
$\mathcal{P}_{\Omega}(\bY)=\bm{\Omega}\odot(\bX+\bO+\bV)$
[cf. \eqref{eq:data_model}], and smart meter $n$ has available the $n$-th
row of $\mathcal{P}_{\Omega}(\bY)$.

\begin{figure}[t]
\centering
\includegraphics[width=\linewidth]{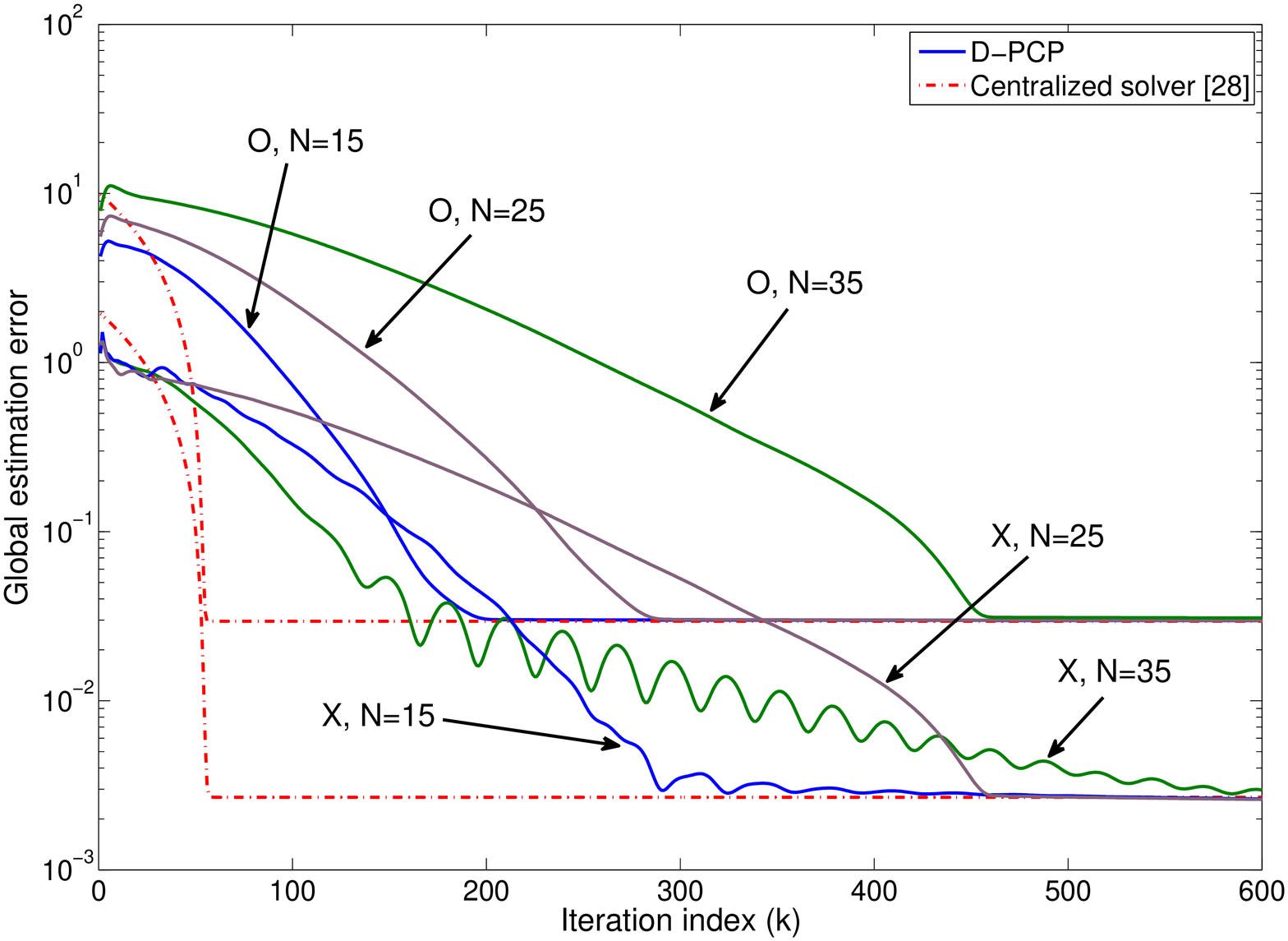}
\caption{Convergence of the D-PCP algorithm for different network sizes. D-PCP 
attains the same estimation error as the centralized solver.}
\label{fig:Fig_2}
\end{figure}

\begin{figure}[t]
\centering
\includegraphics[width=\linewidth]{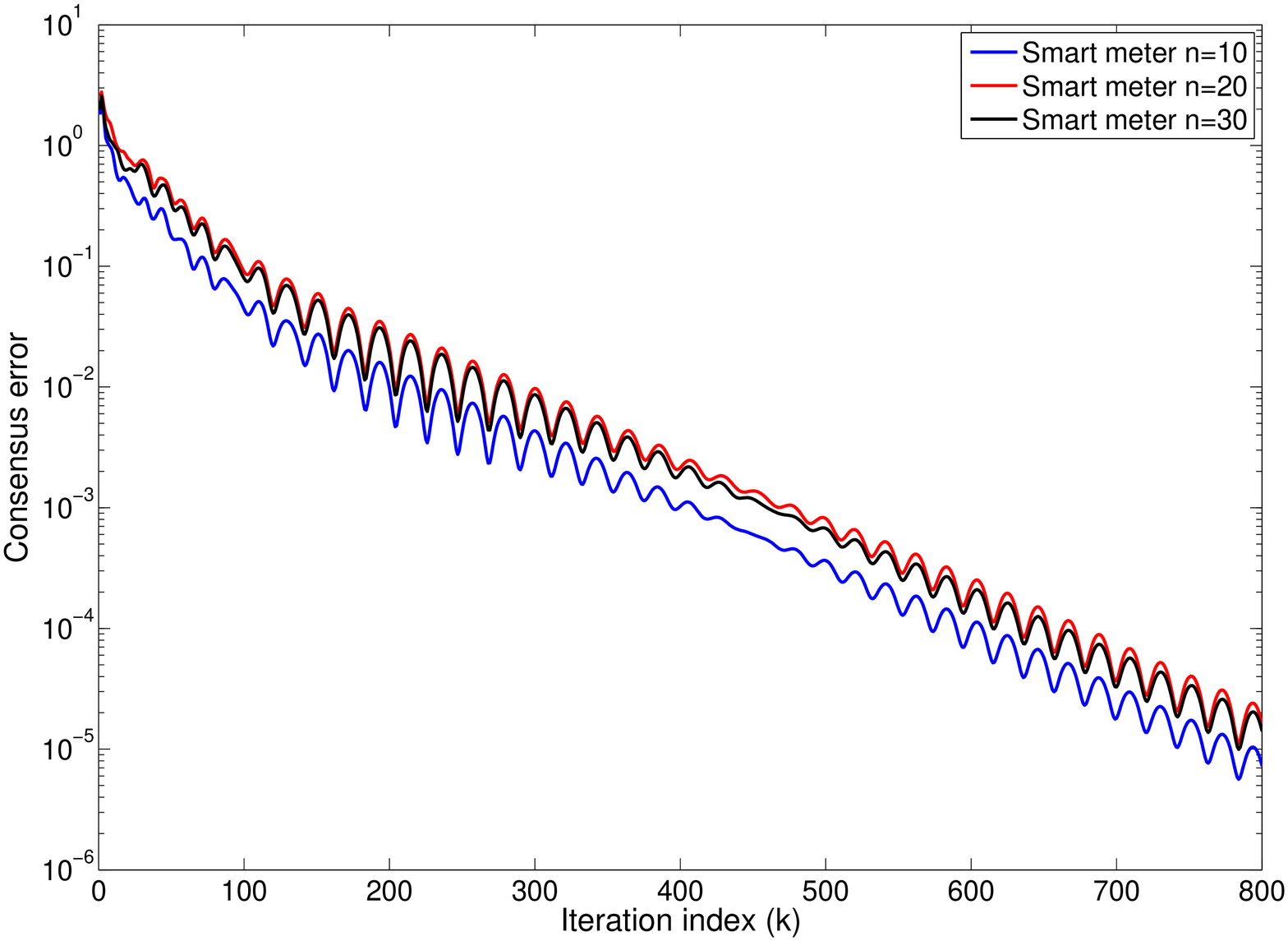}
\caption{Evolution of the consensus error.}
\label{fig:Fig_3}
\vspace{-0.5cm}
\end{figure}

To experimentally corroborate the convergence and optimality 
(as per Proposition \ref{prop:ther_1}) of the D-PCP algorithm, Algorithm \ref{tab:table_2} is
run with $c=1$ and compared with the centralized benchmark
(P1), obtained using the solver in~\cite{rpca_admom}. Parameters $\lambda_1=0.0141$ and
$\lambda_*=0.346$ are chosen as suggested in~\cite{zlwcm10}. For both
schemes, Fig. \ref{fig:Fig_2} shows the evolution of the global estimation errors 
$e_X[k]:=\|\bX[k]-\bX\|_F/\|\bX\|_F$ and $e_O[k]:=\|\bO[k]-\bO\|_F/\|\bO\|_F$.
It is apparent that the D-PCP algorithm converges to the centralized estimator,
and as expected convergence slows down due to the delay associated
with the information flow throughout the network. The test is also repeated
for network sizes of $N=15$ and $35$ devices, to illustrate that the 
time till convergence scales gracefully as the network size increases. Finally,
for $N=35$ and with $\bar{\bQ}[k]:=\sum_{n}\bQ_n[k]/N$, Fig. \ref{fig:Fig_3} depicts the consensus error 
$e_{c,n}[k]:=\|\bQ_n[k]-\bar{\bQ}[k]\|_F/\|\bar{\bQ}[k]\|_F$ for three representative
smart metering devices. In all cases the error decays rapidly to zero, showing
that networkwide agreement is attained on the estimates $\bQ_n[k]$


\subsection{Real load curve data test}
\label{ssec:real_data}

\begin{figure}[t]
\centering
\includegraphics[width=\linewidth]{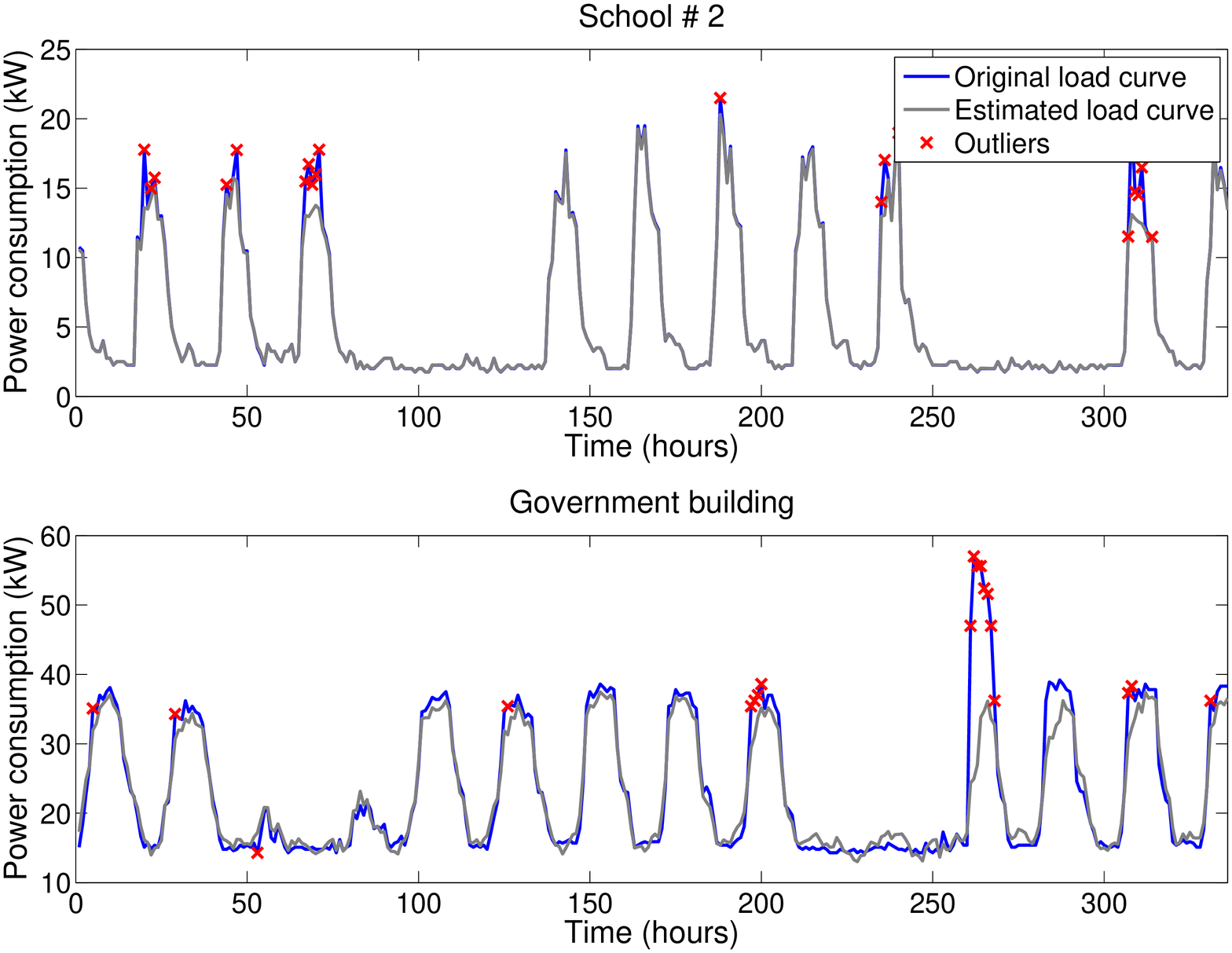}
\caption{School and government building load curve data cleansing.}
\label{fig:Fig_4}
\vspace{-0.5cm}
\end{figure}

\begin{figure}[t]
\centering
\includegraphics[width=\linewidth]{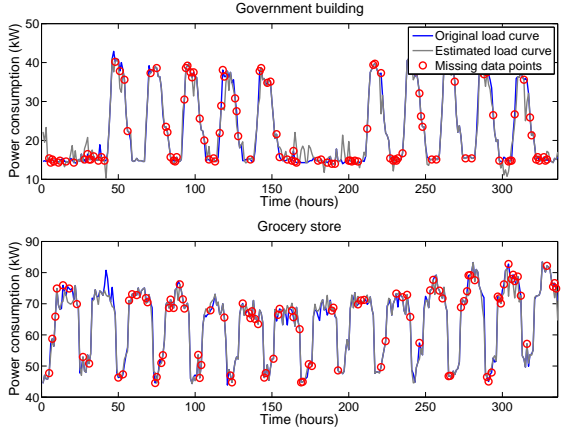}
\caption{Government building and grocery store load curve imputation, when $30\%$ of the
data are missing.}
\label{fig:Fig_5}
\vspace{-0.5cm}
\end{figure}

Here, the D-PCP algorithm is tested 
on real load curve data. The dataset consists of power consumption measurements (in kW) 
for a government building, a grocery store, and three schools ($N=5$)
collected every fifteen minutes during a period of 
more than five years, ranging from July 2005 to October 2010. Data is 
downsampled by a factor of four, to yield one measurement per hour. For the present experiment, 
only a subset of the whole data is utilized for concreteness, where $T=336$ 
was chosen corresponding to $336$ hour periods. For the government building
case, a snapshot of the available load curve data spanning the studied
two-week period is shown in blue e.g., in Fig. \ref{fig:Fig_4} (bottom). Weekday activity 
patterns can be clearly discerned from those corresponding to weekends, as expected for
most government buildings; but different, e.g., for the load profile of the
grocery store in Fig. \ref{fig:Fig_5} (bottom). 

To run the D-PCP algorithm, an underlying communication graph was generated
as in Section \ref{ssec:sim_data}. A randomly chosen subset of $30\%$ of the measurements
was removed to model missing data.
For one of the schools and the government building data, Fig. \ref{fig:Fig_2} depicts the 
cleansed load curves that closely follow the measurements, but are smooth 
enough to avoid overfitting the abnormal energy peaks on the so-termed ``building
operational shoulders.'' Indeed, these peaks are in most cases identified as outliers.
The effectiveness in terms of imputation of missing data is illustrated in Fig. \ref{fig:Fig_3}
(identified outliers are not shown here); note how the cleansed (gray) load curve
goes through the (red) missing data points. The relative error in predicting missing
data is around $6\%$, and degrades to $8\%$ when the amount of missing data
increases to $50\%$.

\section{Conclusion}\label{sec:conc}
A novel robust load curve cleansing and imputation method
is developed in this paper, rooted at the crossroads of 
sparsity-cognizant statistical inference, low-rank matrix completion, 
and large-scale distributed optimization. The adopted PCP estimator
jointly leverages the low-intrinsic dimensionality of spatiotemporal load
profiles, and the sparse (that is, sporadic) nature of outlying measurements.
A separable reformulation of PCP is shown to be efficiently minimized using the ADMM,
and gives rise to fully-decentralized iterations which can be run
by a network of smart-metering devices. Comprehensive tests with synthetic
and real load curve data demonstrate the effectiveness of the novel load
cleansing and imputation approach, and corroborate the convergence
and global optimality of the D-PCP algorithm.

An interesting future direction is to devise real-time cleansing
and imputation algorithms capable of processing load curve data acquired sequentially in
time. Online adaptive algorithms enable tracking of ``bad data'' in nonstationary
environments, typically arising due to to e.g., network topology changes and missing data. 
In addition, it is of interest to rigorously establish
convergence of the D-PCP algorithm. Such results
could significantly broaden the applicability of ADMM for large-scale
optimization over networks, even in the presence of non-convex but highly 
structured and separable cost functions.



\appendix[Algorithmic Construction]
The goal is to show that [S1]-[S4] can be simplified to
the iterations tabulated under Algorithm \ref{tab:table_2}. 
Focusing first on [S3], \eqref{S3_ADMOM} decomposes into $N$
ridge-regression sub-problems
\begin{equation*}
\bl_n[k+1]=\arg\min_{\bl}\left\{\left\|\bby_{n} 
-\bQ_{n}[k+1]\bl-\bbo_n[k]\right\|_2^2+\lambda_*\|\bl\|_2^2\right\}\nonumber
\end{equation*}
which admit the closed-form solutions shown in Algorithm \ref{tab:table_2}.

Moving on to [S4], from the decomposable structure of the augmented Lagrangian
[cf. \eqref{augLagr}] \eqref{S4_ADMOM} decouples into per-node
scalar Lasso subtasks (note that $\bQ_n:=[\bq_{n,1},\ldots,\bq_{n,T}]'$)
\begin{align}
o_{n,t}[k+1]={}&\arg\min_{\bbo}\left\{\frac{1}{2}\left(y_{n,t} 
-\bq_{n,t}'[k+1]\bl_n[k+1]-o\right)^2\right.\nonumber\\
&\hspace{1cm}\left.\phantom{\frac{1}{2}}+ \lambda_1|o|_1\right\}, \quad t=1,\ldots,T\nonumber\\
={}&  \mathcal{S}_{\lambda_1}(y_{n,t} -\bq_{n,t}'[k+1]\bl_n[k+1]),~t=1,\ldots,T\nonumber
\end{align}
and 
$\sum_{n=1}^{N}|\cJ_n|$ additional unconstrained QPs
\begin{align}\label{gammalocalproblem}
\nonumber\hspace{-0.2cm}\bar{\bF}_n^{m}[k+1]=\tilde{\bF}_n^{m}[k+1]=\mbox{arg}\:\min_{\bar{\bF}_n^{m}}
&\left\{-\langle\bar{\bM}_{n}^{m}[k]+\tilde{\bM}_{n}^{m}[k],\bar{\bF}_n^{m}\rangle\phantom{\frac{c}{2}}\right.\\
&\hspace{-5cm}\left.+ \frac{c}{2}\left[\|\bQ_{n}[k+1]-\bar{\bF}_n^{m}\|_F^2+
\|\bQ_{m}[k+1]-\bar{\bF}_n^{m}\|_F^2\right]\right\}
\end{align}
which admit the closed-form solutions
\begin{align}\label{gammalocalrecursion}
\nonumber\bar{\bF}_n^{m}[k+1]=\tilde{\bF}_n^{m}[k+1]=&\frac{1}{2c}
(\bar{\bM}_{n}^{m}[k]+\tilde{\bM}_{n}^{m}[k])\\
&\hspace{-1cm}+\frac{1}{2}\left(\bQ_{n}[k+1]+\bQ_{m}[k+1]\right).
\end{align}
Note that in formulating \eqref{gammalocalproblem},
$\tilde{\bF}_n^{m}$ was eliminated using the constraints
$\nonumber\bar{\bF}_n^{m}=\tilde{\bF}_n^{m}$ defining $\mathcal{C}_F$. Using
\eqref{gammalocalrecursion} to eliminate
$\bar{\bF}_n^{m}[k]$ and
$\tilde{\bF}_n^{m}[k]$ from \eqref{eq:multi_barC} and
\eqref{eq:multi_tildeC} respectively, a simple induction argument
establishes that if the initial Lagrange multipliers obey
$\bar\bM_{n}^{m}[0]=-\tilde{\bM}_{n}^{m}[0]=\mathbf{0}$,
then $\bar\bM_{n}^{m}[k]=-\tilde{\bM}_{n}^{m}[k]$ for all
$k\geq 0$, where $n\in\calN$ and $m\in\calJ_n$. The set
$\{\tilde{\bM}_{n}^{m}\}$ of multipliers has been shown
redundant, and \eqref{gammalocalrecursion} readily simplifies
to
\begin{equation}
\bar{\bF}_n^{m}[k+1]=\tilde{\bF}_n^{m}[k+1]=
\frac{1}{2}\left(\bQ_{n}[k+1]+\bQ_{m}[k+1]\right).
\label{gammalocalrecursionsimple}
\end{equation}
It then follows that
$\bar{\bF}_n^{m}[k]=\bar{\bF}^n_{m}[k]$ for
all $k\geq0$, an identity that will be used later on. By
plugging \eqref{gammalocalrecursionsimple} in \eqref{eq:multi_barC}, the
(non-redundant) multiplier updates become
\begin{equation}\label{vupdate}
\bar\bM_{n}^{m}[k]=\bar\bM_{n}^{m}[k-1]+\frac{c}{2}[\bQ_{n}[k]-\bQ_{m}[k]],
{\quad}n\in\calN,{\:}m\in\calJ_n.
\end{equation}
If $\bar\bM_{n}^{m}[0]=-\bar\bM_{m}^{n}[0]=\mathbf{0}$,
then the structure of \eqref{vupdate} reveals that
$\bar\bM_{n}^{m}[k]=-\bar\bM_{m}^{n}[k]$ for all $k\geq
0$, where $n\in\calN$ and $m\in\calJ_n$.

The minimization \eqref{S4_ADMOM} in [S4] also decouples in $N$
simpler sub-problems, namely
\begin{align}
\nonumber\bQ_n[k+1]=\arg\min_{\bQ}&
\left\{\frac{1}{2}\|\bm\Omega_n(\mathbf{y}_n 
-\bQ\bl_n[k]-\mathbf{o}_n[k])\|_2^2
\right.\nonumber\\
&\hspace{-1.5cm}+\frac{\lambda_*}{2N}\|\bQ\|_F^2+\sum_{m\in\calJ_{n}}
\langle\bar{\bM}_{n}^{m}[k]+\tilde{\bM}_{m}^{n}[k],\bQ\rangle
\nonumber\\
&\hspace{-1.5cm}\left.+\frac{c}{2}\sum_{m\in\calJ_{n}}
\left(\|\bQ-\bar{\bF}_{n}^{m}[k]\|_F^{2}+
\|\bQ-\tilde{\bF}_{m}^{n}[k]\|_F^2\right)\right\}\nonumber\\
\nonumber=\arg\min_{\bQ}&\left\{
\frac{1}{2}\|\bm\Omega_n(\mathbf{y}_n 
-\bQ\bl_n[k]-\mathbf{o}_n[k])\|_2^2\right.\\
&\hspace{-3cm}\left.+\frac{\lambda_*}{2N}\|\bQ\|_F^2+\langle\bS_n[k],\bQ\rangle
+c\sum_{m\in\calJ_{n}}
\left\|\bQ-\frac{\bQ_n[k]+\bQ_m[k]}{2}\right\|_F^{2}\right\}\label{eq:Q_update}
\end{align}
where in deriving \eqref{eq:Q_update} it was used that: i)
$\bar\bM_{n}^{m}[k]=\tilde{\bM}^{n}_{m}[k]$ which follows
from the identities
$\bar\bM_{n}^{m}[k]=-\tilde{\bM}_{n}^{m}[k]$ and
$\bar\bM_{n}^{m}[k]=-\bar\bM_{m}^{n}[k]$ established
earlier; ii) the definition
$\bS_n[k]:=2\sum_{m\in\calJ_n}\bar\bM_{n}^{m}[k]$; and
iii) the identity
$\bar{\bF}_{n}^{m}[k]=\tilde{\bF}^{n}_{m}[k]$
which allows one to merge the identical quadratic penalty terms and
eliminate both $\bar{\bF}_{n}^{m}[k]$ and
$\tilde{\bF}^{n}_{m}[k]$ using
\eqref{gammalocalrecursionsimple}. Problem (14)
is again an unconstrained QP, which is readily solved
in closed form by e.g., vectorizing $\bQ$ and examining the
first-order condition for optimality. 

Finally, note that upon scaling by two the
recursions \eqref{vupdate} and summing them over
$m\in\calJ_n$, the update recursion for $\bS_n[k]$ in
Algorithm \ref{tab:table_2} follows readily.\hfill$\blacksquare$

\section*{Acknowledgment}

The authors would like to thank NorthWrite
Energy Group and Prof. Vladimir Cherkassky 
(Dept. of ECE, University of Minnesota) for providing
the data analysed in Section \ref{ssec:real_data}.

\ifCLASSOPTIONcaptionsoff
  \newpage
\fi



%

%
\vspace{-0.5cm}
\begin{biography}[{\includegraphics[width=1in,height=1.25in,clip,keepaspectratio]{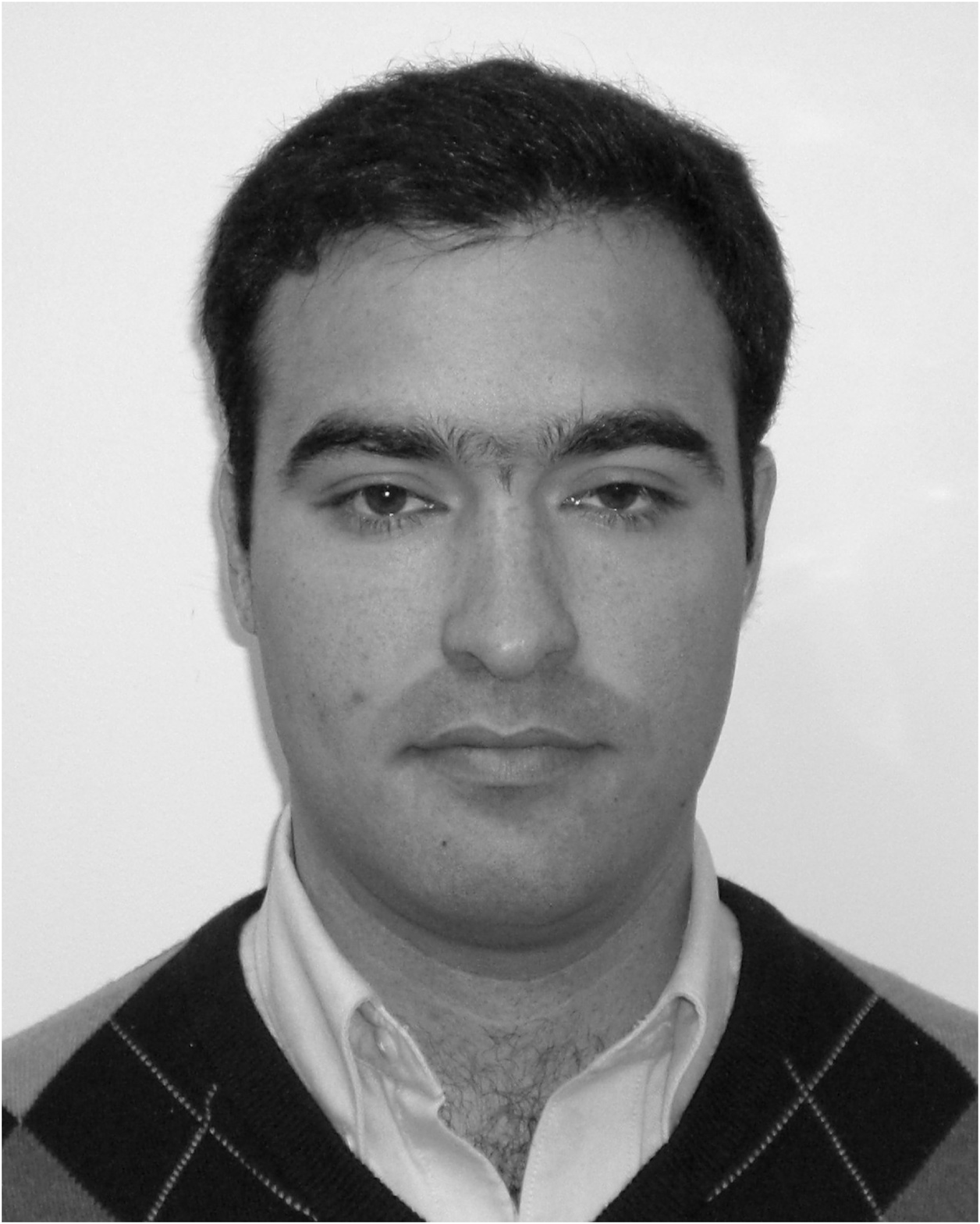}}]
{Gonzalo Mateos (M'12)} 
received his B.Sc. degree in Electrical Engineering 
from Universidad de la Rep\'{u}blica (UdelaR), Montevideo, Uruguay in 2005 
and the M.Sc. and Ph.D. degrees in Electrical and Computer
Engineering from the Univ. of Minnesota, Minneapolis, in 2009 and 2011. 
Since 2012, he has been a post doctoral
associate with the Dept. of Electrical and Computer Engineering
and the Digital Technology Center, Univ. of Minnesota.

From 2003 to 2006, he was a teaching assistant with the Dept. of Electrical Engineering,  
UdelaR. From 2004 to 2006, he worked as a Systems Engineer at 
Asea Brown Boveri (ABB), Uruguay. His research interests lie in 
the areas of communication theory, signal processing and networking. 
His current research focuses on distributed optimization, sparsity-cognizant
signal processing, and statistical learning for cartography of
cognitive networks including the smart power grid.
\end{biography}
\vspace{-0.5cm}
\begin{biography}[{\includegraphics[width=1in,height=1.25in,clip,keepaspectratio]{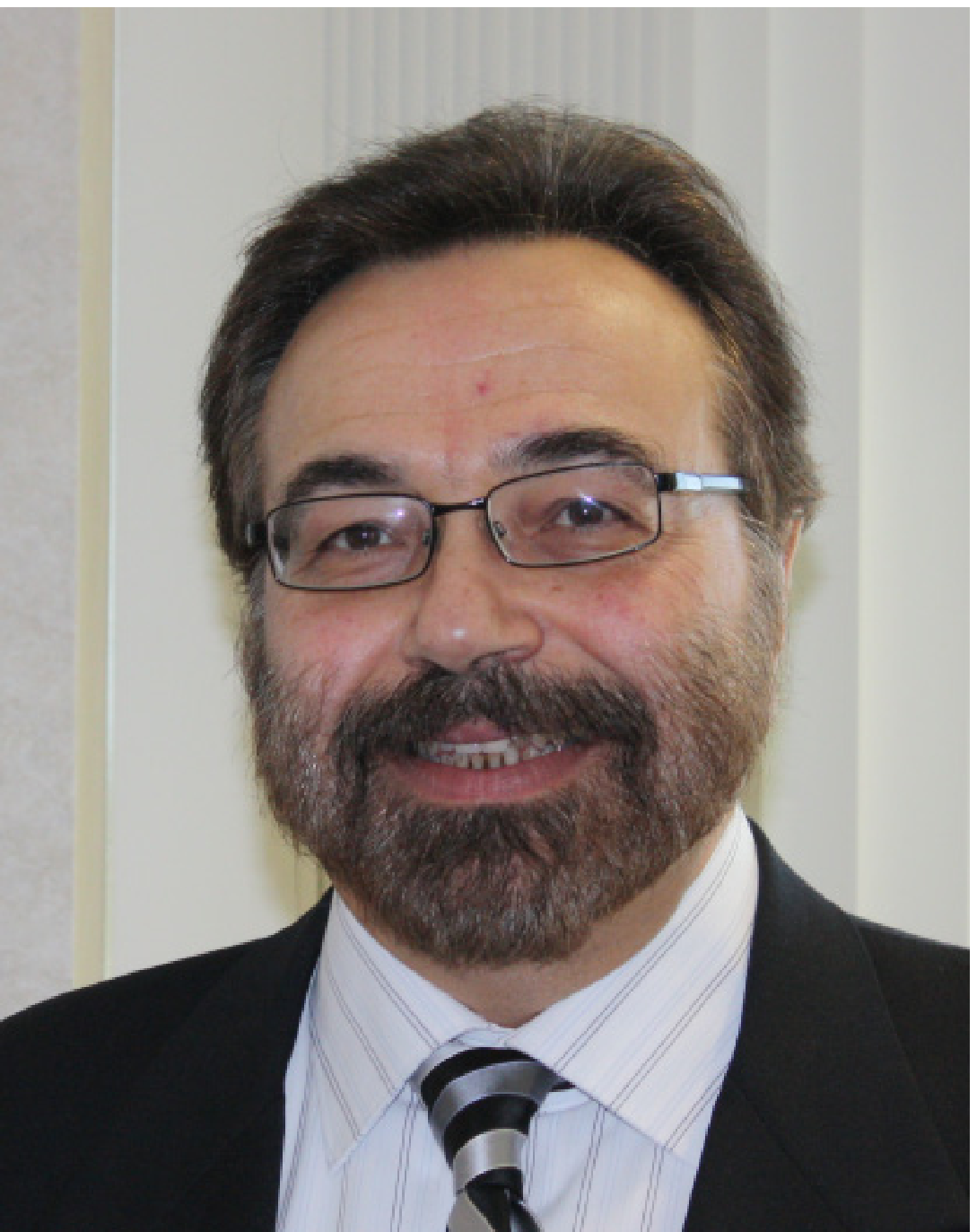}}]
{Georgios B.  Giannakis (Fellow'97)} 
received his Diploma in Electrical Engr. from the Ntl. Tech. Univ. of Athens, 
Greece, 1981. From 1982 to 1986 he was with the Univ. of Southern California 
(USC), where he received his MSc. in Electrical Engineering, 1983, MSc. in 
Mathematics, 1986, and Ph.D. in Electrical Engr., 1986. Since 1999 he has 
been a professor with the Univ. of Minnesota, where he now holds an ADC 
Chair in Wireless Telecommunications in the ECE Department, and serves as 
director of the Digital Technology Center.

His general interests span the areas of communications, networking and 
statistical signal processing - subjects on which he has published more 
than 340 journal papers, 560 conference papers, 20 book chapters, two edited 
books and two research monographs. Current research focuses on compressive 
sensing, cognitive radios, cross-layer designs, wireless sensors, social and 
power grid networks. He is the (co-) inventor of 21 patents issued, and the 
(co-) recipient of 8 best paper awards from the IEEE Signal Processing (SP) 
and Communications Societies, including the G. Marconi Prize Paper Award in 
Wireless Communications. He also received Technical Achievement Awards from 
the SP Society (2000), from EURASIP (2005), a Young Faculty Teaching Award, 
and the G. W. Taylor Award for Distinguished Research from the University 
of Minnesota. He is a Fellow of EURASIP, and has served the IEEE in a 
number of posts, including that of a Distinguished Lecturer for the IEEE-SP Society.
\end{biography}




\end{document}